\documentclass[reqno, 12pt]{amsart}
\usepackage[letterpaper,hmargin=1in,vmargin=1in]{geometry}
\usepackage{amsmath, amssymb, amsthm, verbatim, bbm, url} 

\newcommand \R {\mathbb{R}}
\newcommand \C {\mathbb{C}}
\newcommand \Z {\mathbb{Z}}

\newcommand \Oh {\mathcal{O}}

\newcommand \D {\partial}

\newcommand \Def {\stackrel{\textrm{def}}=}

\DeclareMathOperator \supp {supp}

\DeclareMathOperator \Hess {Hess}

\DeclareMathOperator \Tr {Tr}

\newtheorem*{thm}{Theorem}
\newtheorem*{cor}{Corollary}

\theoremstyle{definition}

\numberwithin{equation}{section}
\numberwithin{lem}{section}
\numberwithin{Defn}{section}

\title
{Spectral uniqueness of radial semiclassical Schr\"odinger operators}

\author[Kiril Datchev]
{Kiril Datchev}
\address{Department of Mathematics, Massachusetts Institute of Technology, Cambridge, MA
02139.}
\email{datchev@math.mit.edu}
\author[Hamid Hezari]
{Hamid Hezari}
\address{Department of Mathematics, Massachusetts Institute of Technology, Cambridge, MA 02139.}
\email{hezari@math.mit.edu}
\author[Ivan Ventura]
{Ivan Ventura}
\address{Department of Mathematics, University of California, Berkeley, CA 94720.}
\email{iventura@math.berkeley.edu}

\thanks{The first author is partially supported by a National Science Foundation postdoctoral fellowship, and the second author is partially supported by the National Science Foundation under
grant DMS-0969745. The authors are grateful for the hospitality of the Mathematical Sciences Research Institute, where part of this research was carried out.}
\date{February 11, 2011}

\begin{document}
\begin{abstract}
We prove that the spectrum of an $n$-dimensional semiclassical radial Schr\"odinger operator determines the potential within a large class of potentials  for which we assume no symmetry or analyticity. Our proof is based on the first two semiclassical trace invariants and on the isoperimetric inequality.
\end{abstract}
\maketitle

\section{Introduction}

In this paper we study the problem of determining a potential $V$ in the class
\begin{equation}\label{e:class} \{ V \in C^\infty(\R^n):  -h^2\Delta + V \textrm{ is selfadjoint on a domain containing }C_0^\infty(\R^n)\} \end{equation}
from its eigenvalues. Note that this class includes the set of potentials $V$ which are bounded from below, because in this case the Friedrichs extension satisfies the condition.

We prove that radial, monotonic potentials are spectrally determined within this class.

Let $n \ge 2$. Let $V_0 \in C^\infty(\R^n)$ be radial near $0$ in the sense that $V_0(x) = R(|x|)$ for $|x| \le R_0$ for some $R_0>0$, and suppose $R$ satisfies $R(0)=0$ and $R'(r) > 0$ for $r \in (0,R_0)$. Let $\lambda_0 = R(R_0)$, and suppose $V_0 > \lambda_0$ for $|x| > R_0$.

\begin{thm}
If $V$ is in the class \eqref{e:class} and  the spectrum of $-h^2\Delta + V$ agrees with the spectrum of $-h^2\Delta + V_0$ in the interval $(-\infty,\lambda_0)$ up to order $o(h^2)$ for $h \in (0,h_0]$ for some $h_0>0$, then there exists $x_0 \in \R^n$ such that $V(x-x_0)=V_0(x)$ for $|x| < R_0$.
\end{thm}

In particular, the harmonic oscillator is spectrally determined.

\begin{cor}
If $V$ is in the class \eqref{e:class} and if the spectrum of $-h^2\Delta + V$ agrees with the spectrum of $-h^2\Delta + |x|^2$ up to order $o(h^2)$ for $h \in (0,h_0]$ for some $h_0>0$, then there exists $x_0 \in \R^n$ such that $V(x-x_0)=|x|^2$.
\end{cor}

The oldest inverse spectral result is that balls in $\R^n$ are determined by their Dirichlet or Neumann eigenvalues within the class of bounded domains with smooth boundary. This follows from the isoperimetric inequality and the fact that the first two heat invariants, i.e. the first two terms of the asymptotic expansion of $\Tr e^{t\Delta}$ as $t \to 0^+$, give the volume of the domain and the area of its boundary. All other positive results for spectral uniqueness of domains rely on analyticity or symmetry assumptions on the domain and on the class of domains under consideration. In \cite{Z09}, Zelditch shows that any analytic domain in $\R^2$ with one reflection symmetry is determined within the class of all such domains, and in \cite{HZ10a}, the second author and Zelditch generalize this result to higher dimensions. In \cite{HZ10}, the second author and Zelditch show that an isospectral deformation of an ellipse within the class of smooth domains with two reflection symmetries is infinitesimally trivial.

In general it is impossible to determine a potential $V$ from the spectrum of the nonsemiclassical Schr\"odinger operator $-\Delta + V$. For example, in \cite{MT81}, McKean and Trubowitz find an infinite dimensional family of potentials in $C^\infty(\R)$ which are isospectral with the harmonic oscillator $V(x) = x^2$. However, if we have information about the spectrum of the semiclassical Schr\"odinger operator  $-h^2\Delta + V$ for all $h$ in $(0,h_0]$, we can say more about the potential. In fact, in the Theorem above, it is sufficient for example to know the spectrum up to $o(h^2)$ for $h \in \{h_j\}_{j=1}^\infty$ for a sequence $h_j \to 0$. In a physical interpretation of the problem, varying $h$ corresponds to varying a parameter of the system, such as a mass or coupling constant.

The semiclassical Schr\"odinger operator is, in the following sense, the natural spectral analogue of the Laplacian on a domain or manifold. A semiclassical parameter $h$ can be introduced into any spectral problem about the Laplacian on a domain or manifold by scaling the metric by $h^{-2}$, and hence producing a factor of $h^2$ in front of the Laplacian. Because of the homogeneity of the symbol of the Laplacian in the fibers of the cotangent bundle, the spectrum of the semiclassical Laplacian for any value of $h$ can be deduced from the spectrum at $h=1$:
\begin{equation}\label{e:intro} -\Delta u_j = \lambda_j u_j \quad\Longleftrightarrow\quad -h^2 \Delta u_j = E_j(h) u_j,  \textrm{ for \textit{all} $h>0$, where $E_j(h) = h^2 \lambda_j$}.\end{equation}
The symbol of a Schr\"odinger operator is not homogeneous in the fibers, and so in this setting there is no analogue of \eqref{e:intro}. This partially accounts for the large family of isospectral potentials in the example of McKean and Trubowitz mentioned above: the nonsemiclassical spectrum alone provides much less information in this setting.

In \cite{ISZ}, Iantchenko, Sj\"ostrand and Zworski study spectral inverse problems of semiclassical Schr\"odinger operators, recovering in particular Birkhoff normal forms of closed orbits. In \cite{gu07}, Guillemin and Uribe show that any real analytic potential $V$, symmetric with respect to all coordinate axes and with a unique global minimum, is determined within the class of all such potentials by the spectrum at the bottom of the well of its semiclassical Schr\"odinger operator. In \cite{h09}, the second author removes one of the symmetry assumptions, giving a result with no symmetry assumptions in dimension 1, but keeps the assumption of analyticity. In \cite{cg}, Colin de Verdi\`ere and Guillemin give a different proof of the one dimensional result. In \cite{c}, Colin de Verdi\`ere removes both the symmetry and analyticity assumptions in the one dimensional case, but adds a genericity assumption. In \cite{gw10} (see also \cite[\S 10.6]{gs}) Guillemin and Wang  prove a one dimensional version of the result in the present paper: using the two trace invariants
\begin{equation}\label{e:gsinv}\int_{\{|\xi|^2 + V(x) < \lambda\}} dxd\xi, \qquad \int_{\{|\xi|^2 + V(x) < \lambda\}} |\nabla V(x)|^2dxd\xi,\end{equation}
they show that an even function (or a suitable noneven function) is determined by its spectrum within the class of functions monotonic away from $0$. To the authors' knowledge, the result in the present paper is the first positive result for a semiclassical Schr\"odinger operator in higher dimensions with no analyticity, symmetry or genericity assumption on the class within which the potential is spectrally determined.

We mention also some related results for determining an obstacle from its resonances. In \cite{hz99}, Hassell and Zworski show that the $3$-ball as an obstacle is determined by its Dirichlet resonances. In \cite{chr}, Christiansen extends the result to $n$-balls and Neumann resonances.

Our proof is based on the isoperimetric inequality and on the first two terms of a semiclassical trace formula due to Helffer and Robert \cite{hr}. In this it is very similar to the proof of spectral uniqueness of balls in $\R^n$, which uses the isoperimetric inequality and the first two terms of the heat trace. Because we use only the first two terms of the expansion, we only need to know the spectrum up to $o(h^2)$.

We are grateful to Maciej Zworski for his encouragement, and for several interesting discussions about this paper and about the motivation for the inverse spectral problem for semiclassical Schr\"odinger operators.

\section{Review of semiclassical trace formulas}

For $V \in C^\infty(\R^n;[0,\infty))$, let
\[P \Def -h^2 \Delta + V.\]
If $P$ has discrete spectrum in $[0,\lambda_0)$, then for each $h>0$, $f(P)$ (defined by the functional calculus for selfadjoint operators) has finite rank for every $f \in C_0^\infty((-\infty,\lambda_0))$, and the trace of $f(P)$ is determined by the eigenvalues of $P$. In this section we review three methods of proof for the trace formula
\begin{equation}\label{e:trace}
\begin{split}\Tr(&f(P)) = \\&\frac 1 {(2\pi h)^n} \left(\int_{\R^{2n}} f(|\xi|^2 + V)dxd\xi + \frac{h^2}{12} \int_{\R^{2n}} |\nabla V|^2 f^{(3)}(|\xi|^2 + V)dxd\xi + \Oh(h^4)\right).\end{split}
\end{equation}
The two invariants in \eqref{e:gsinv} are obtained from this formula by approximating the characteristic function of an interval in $C_0^\infty(\R)$.

The trace formula \eqref{e:trace} comes from an asymptotic expansion of $\sigma(f(P))$, the full semiclassical symbol of $f(P)$ (here the Weyl symbol, or the standard symbol, or another symbol may be used), of the form
\begin{equation}\label{e:symbolexp}\sigma(f(P)) = \sigma_0(f(P)) + h \sigma_1(f(P)) + h^2\sigma_2(f(P)) + \cdots\end{equation}
combined with a version of Lidskii's formula for the trace of an operator as the integral along the diagonal of its integral kernel, as in for example \cite[Theorem 9.6]{ds}:
\[\Tr(f(P)) = \frac 1 {(2\pi h)^n} \left(\int_{\R^{2n}} \left[\sigma_0(f(P)) + h \sigma_1(f(P)) + h^2\sigma_2(f(P)) + \cdots\right]dxd\xi\right).\]

The terms of the expansion \eqref{e:symbolexp} were first computed by Helffer and Robert \cite[Proposition 5.3]{hr} using the Mellin inversion formula
\[f(P) = \frac 1{2\pi i} \int_{\sigma-i\infty}^{\sigma+i\infty} F(s) P^{-s} ds,\]
where $F$ is the Mellin transform of $f$, given by $F(s) = \int_0^\infty t^{s-1}f(t)dt$. In \cite[Theorem 8.7]{ds} Dimassi and Sj\"ostrand obtain an expansion using the Helffer-Sj\"ostrand formula (\cite{hs}, see also for example \cite[(8.3)]{ds})
\[f(P) = - \frac 1 \pi \int_\C \overline\D \tilde f(z)(z-P)^{-1}dxdy,\]
where $\tilde f$ is an almost analytic extension of $f$. More recently, in \cite[\S 10]{gs}, Guillemin and Sternberg use the Fourier inversion formula
\[f(P) =  \frac 1 {2\pi} \int_\R \hat f(t) e^{itP}dt\]
to obtain the same expansion. This last method is a semiclassical version of Taylor's construction in \cite[\S 12.1]{tay}. In each case this reduces the problem to computing the asymptotics of the symbol of one of $P^{-s}$, $(z-P)^{-1}$ or $e^{itP}$ (in fact in the first case Helffer and Robert use the contour integral formula $P^{-s} = (2\pi i)^{-1} \int_\gamma z^{-s} (z - P)^{-1}dz$ for suitable $\gamma$ to reduce the problem to studying $(z-P)^{-1}$). Using one of these methods, and either the standard or Weyl symbol, we obtain an expansion of the form \eqref{e:symbolexp}. Regardless of the choice of symbol, we have $\sigma_0(f(P)) = f (\sigma_0(P))$, i.e. the principal symbol of $f(P)$ is $f$ of the principal symbol of $P$. Using the formulas for $\sigma_0, \sigma_1,\sigma_2,\sigma_3$ derived in \cite[\S 10.5]{gs} gives \eqref{e:trace}.

\section{Proof of Theorem}\label{s:main}
The spectrum of $-h^2 \Delta + V_0$ is contained in $(0,\infty)$ and is discrete in $[0,\lambda_0)$. Since the spectrum of $-h^2 \Delta + V$ is the same as the spectrum of $-h^2 \Delta + V_0$, up to order $o(h^2)$ on $(-\infty,\lambda_0)$, we conclude that $0 \le -h^2\Delta + V  + o(h^2)$. It follows that $V \ge 0$. Indeed, otherwise we could take a function $u \in C_0^\infty(\R^n)$ with $\supp u \subset \{V < 0\}$, giving
\[0 \le \int (-h^2 \Delta + V)|u|^2dx + o(h^2)\|u\|^2_{L^2} = \int_{\supp u} V |u|^2 + \Oh(h^2)\|u\|_{H^1} < 0,\]
for $h$ sufficiently small, which is a contradiction.

Because the spectrum of $-h^2 \Delta + V$ is known on $(-\infty,\lambda_0)$ up to order $o(h^2)$ we obtain from \eqref{e:trace} the two trace invariants in \eqref{e:gsinv}
for each $\lambda \in (0,\lambda_0)$. By integrating in the $\xi$ variable, we rewrite these invariants as follows:
\begin{equation}\label{e:same}\int_{\{V(x) < \lambda\}} (\lambda - V)^{n/2} dx, \qquad \int_{\{V(x) < \lambda\}} |\nabla V(x)|^2 (\lambda - V)^{n/2} dx.\end{equation}

From the first invariant of \eqref{e:same} it follows that
\begin{equation}\label{e:critset} \{x\colon V(x) < \lambda, \nabla V(x) = 0\} \textrm{ has Lebesgue measure zero.} \end{equation}
We will prove this in \S\ref{s:critset}.

Using \eqref{e:critset} and the coarea formula we rewrite the invariants in \eqref{e:same} as
\[\int_0^\lambda \left(\int_{\{V=s, \nabla V \ne 0\}} \frac{(\lambda - V)^{n/2}}{|\nabla V|} dS\right)ds, \qquad \int_0^\lambda \left(\int_{\{V=s\}} {|\nabla V|}{(\lambda - V)^{n/2}} dS\right)ds.\]

Using the fact that $V=s$ in the inner integrand, the factor of $(\lambda - V)^{n/2} = (\lambda - s)^{n/2}$ can be taken out of the surface integral, leaving
\begin{equation}\label{e:abel}\int_0^\lambda (\lambda - s)^{n/2} I_1(s)ds, \qquad \int_0^\lambda (\lambda - s)^{n/2} I_2(s)ds,\end{equation}
where
\begin{equation}\label{e:surfaceinvariants}I_1(s) = \int_{\{V=s, \nabla V \ne 0\}} \frac 1 {|\nabla V|} dS, \qquad I_2(s) = \int_{\{V=s\}} |\nabla V| dS. \end{equation}

We denote the integrals \eqref{e:abel} by $A_{1 + n/2}(I_1)(\lambda)$ and $A_{1+ n/2}(I_1)(\lambda)$. These are Abel fractional integrals of $I_1$ and $I_2$ (see for example \cite[\S5.2]{Zelditch} and \cite[(10.45)]{gs}), and they can be inverted by applying $A_{1 + n/2}$, using the formula 
\begin{equation}\label{e:abel2}\frac 1 {\Gamma(\alpha)}A_\alpha\circ \frac 1 {\Gamma(\beta)}A_\beta = \frac 1 {\Gamma(\alpha + \beta)}A_{\alpha+\beta},\end{equation}
and differentiating $n + 1$ times. From this we conclude that the functions $I_1$ and $I_2$ in \eqref{e:surfaceinvariants} are spectral invariants for every $s \in (0,\lambda_0)$.

Integrating $I_1$ and using the coarea formula again as well as \eqref{e:critset} we find that the volumes of the sets $\{V < s\}$ are spectral invariants for each $s \in (0,\lambda_0)$:
\begin{equation}\label{e:volumeinvariants}\int_0^s I_1(s')ds' =\int_0^s \int_{\{V=s', \nabla V \ne 0\}} \frac 1 {|\nabla V|} dSds'  = \int_{\{V< s\}} 1 dx.\end{equation}

From Cauchy-Schwarz and the fact that $I_1$ and $I_2$ are spectral invariants we obtain
\begin{equation}\label{e:cs}\left(\int_{\{V=s\}}1 dS\right)^2 \le \int_{\{V=s\}} \frac 1 {|\nabla V|} dS\int_{\{V=s\}} |\nabla V|dS = \int_{\{R=s\}} \frac 1 {R'} dS \int_{\{R=s\}} R'dS,\end{equation}
when $s$ is not a critical value of $V$, i.e. by Sard's theorem for almost every $s \in (0,\lambda_0)$. On the other hand, using the invariants obtained in \eqref{e:volumeinvariants} and the fact that the sets $\{R<s\}$ are balls, by the isoperimetric inequality (see for example \cite[p1188]{Osserman}) we find
\[ \int_{\{R=s\}} 1 dS \le \int_{\{V=s\}}1 dS.\]
However,
\[\left(\int_{\{R=s\}} 1 dS\right)^2 =  \int_{\{R=s\}} \frac 1 {R'} dS \int_{\{R=s\}} R'dS,\]
because $1/R'$ and $R'$ are constant on $\{R=s\}$. Consequently
\[\int_{\{R=s\}} 1 dS = \int_{\{V=s\}}1 dS,\]
and so $\{V=s\}$ is a sphere for almost every $s$, because only spheres extremize the isoperimetric inequality. Moreover,
\[\left(\int_{\{V=s\}}1 dS\right)^2 = \int_{\{V=s\}} \frac 1 {|\nabla V|} dS \int_{\{V=s\}} |\nabla V|dS,\]
and so $|\nabla V|^{-1}$ and $|\nabla V|$ are proportional on the surface $\{V = s\}$ for almost every $s$, again by Cauchy-Schwarz. Using the equation \eqref{e:cs} to determine the constant of proportionality, we find that
\[|\nabla V|^2 = R'(R^{-1}(s))^2 = (R^{-1})'(s)^{-2} \Def F(s)\]
on $\{V =s\}$. In other words
\begin{equation}\label{e:pde}|\nabla V|^2 = F(V),\end{equation}
for all $x \in V^{-1}(s)$ for almost all $s \in (0,\lambda_0)$. However, because $F(V) \ne 0$ when $V \ne 0$, it follows by continuity that this equation holds for all $x \in V^{-1}((0,\lambda_0))$.

To solve this equation, we restrict it to flowlines of $\nabla V$, with initial conditions taken on a fixed level set $\{V = s_0\}$:
\begin{equation}\label{e:ode}\begin{cases}\dot x(t) = \nabla V(x(t)), \\ x(0) \in \{V = s_0\}. \end{cases}\end{equation}

Observe first that this flow takes level sets to level sets. Indeed,
\[\frac d {dt} V(x(t)) = \nabla V(x(t)) \cdot \dot x(t) = |\nabla V(x(t))|^2 = F(V(x(t))),\]
and this equation can be solved:
\[I(V) \Def \int_{s_0}^V \frac {dV'}{F(V')}, \qquad V(x(t)) = I^{-1}(t).\]
Note that $I$ depends on $s_0$ but not on $x(0)$.

We now use this to show that the integral curves $x(t)$ are lines. Differentiating \eqref{e:ode} we find
\[\ddot x(t) = \frac d {dt} \nabla V(x(t)) = (\Hess V)(x(t)) \cdot \nabla V(x(t)).\]
But differentiating \eqref{e:pde} gives
\[2(\Hess V)\cdot \nabla V = F'(V) \nabla V,\]
so that
\[\ddot x(t) = \frac 12 F'(V(x(t)) \nabla V(x(t)) = \frac 12 F'(I^{-1}(t))\dot x(t).\]
Integrating this equation we find the integral curves are lines:
\[\dot x(t) =  \dot x(0) \frac{\sqrt{F(I^{-1}(t))}}{|\dot x(0)|} = \dot x(0) \frac{\sqrt{F(V(x(t)))}}{|\dot x(0)|}.\]
Moreover, these lines are all directed at the center of the sphere $\{V=s_0\}$ and have a speed dependent only on $s_0$ but not on $x(0)$, and consequently all level sets in the flowout of $\{V = s_0\}$ are spheres with the same center. Since $\nabla V \ne 0$ on $\{0 < V < \lambda_0\}$ thanks to \eqref{e:pde}, the integral curves cover this region, and we conclude that, up to a translation, $V$ is radial on $\{|x| \le R_0\}$. Since the volumes \eqref{e:volumeinvariants} are spectral invariants, it follows that $V(x) = R(|x|)$ on $\{|x| \le R_0\}$.

\section{Analysis of the critical set}\label{s:critset}
In this section we prove \eqref{e:critset}. To do this we will use the first invariant of \eqref{e:same} to show that $V_*dx$, the pushforward of Lebesgue measure by $V$, is absolutely continuous with respect to Lebesgue measure on $(-\infty,\lambda_0)$. From this it follows  that $V^{-1}(E)$ is Lebesgue-null whenever $E \subset (-\infty,\lambda_0)$ is Lebesgue-null. Since the set of critical values of $V$ is Lebesgue-null thanks to Sard's theorem, this completes the proof.

To show that $V_*dx$ is absolutely continuous, we write the first invariant of \eqref{e:same} as
\begin{equation}\label{e:push}\int_{\{V(x) < \lambda\}}(\lambda-V)^{n/2}dx = \int_\R(\lambda-s)^{n/2}  \chi_{[0,\lambda)}(s) (V_* dx),\end{equation}
for every $\lambda \in (0, \lambda_0)$, where $\chi_E$ denotes the characteristic function of $E$. The finiteness of the integrals in \eqref{e:push} implies that $V_*dx(-\infty,\lambda)<\infty$ for each $\lambda < \lambda_0$, and hence $V_*dx$ is a finite, regular, Borel measure on each $(-\infty,\lambda)$. Applying the Abel fractional integral transform as in \S\ref{s:main}, we find that
\[\int_\R(\lambda-s)^n  \chi_{[0,\lambda)}(s) (V_* dx)\]
is a spectral invariant for every $\lambda \in (0,\lambda_0)$. For this we used \eqref{e:abel2}, which is proved using only the Fubini-Tonelli theorem followed by an elementary single-variable integral computation, and hence is valid here.

We now show that, for $m \in\{1,\dots,n\}$,
\begin{equation}\label{e:deriv}\frac d {d\lambda} \int_\R(\lambda-s)^m  \chi_{[0,\lambda)}(s) (V_* dx) = m\int_\R(\lambda-s)^{m-1}  \chi_{[0,\lambda)}(s) (V_* dx),\end{equation}
where the derivative is guaranteed to exist since,by induction, these are each spectral invariants. To prove \eqref{e:deriv}, we write
\[\begin{split}
&\int_\R(\lambda + h -s)^m  \chi_{[0,\lambda + h)}(s)(V_*dx)  - \int_\R(\lambda-s)^m  \chi_{[0,\lambda)}(s) (V_* dx) = \\
&\int_\R \left[(\lambda + h -s)^m  - (\lambda -s)^m\right]\chi_{[0,\lambda + h)}(s)(V_*dx)  + \int_\R(\lambda-s)^m \left[\chi_{[0,\lambda+h)}(s) - \chi_{[0,\lambda)}(s)\right] (V_* dx).
\end{split}\] 
For the first term we have
\[\int_\R \left[(\lambda + h -s)^m  - (\lambda -s)^m\right]\chi_{[0,\lambda + h)}(s)(V_*dx) = m h\int_\R(\lambda -s)^{m-1}\chi_{[0,\lambda + h)}(s)(V_*dx) + \Oh(h^2),\]
thanks to the fact that $\int \chi_{[0,\lambda + h)}(s)(V_*dx) <\infty$ for small $h$. Using the monotone convergence theorem we find that
\begin{equation}\label{e:leftright}
\int_\R(\lambda -s)^{m-1}\chi_{[0,\lambda + h)}(s)(V_*dx) \to
\begin{cases}
\int(\lambda -s)^{m-1}\chi_{[0,\lambda]}(s)(V_*dx), \qquad &h \to 0^+,\\
\int(\lambda -s)^{m-1}\chi_{[0,\lambda)}(s)(V_*dx), \qquad &h \to 0^-.
\end{cases}
\end{equation}
Meanwhile for the second term we have,
\[\int_\R(\lambda-s)^m \left[\chi_{[0,\lambda+h)}(s) - \chi_{[0,\lambda)}(s)\right] (V_* dx) \to
\begin{cases}
\int(\lambda -s)^m\chi_{\{\lambda\}}(s)(V_*dx)=0, \qquad &h \to 0^+,\\
0, \qquad &h \to 0^-.
\end{cases}
\]
Since the derivative in \eqref{e:deriv} must exist, the limits from the left and right in \eqref{e:leftright} must be equal, and \eqref{e:deriv} is proved. In particular, we find that
\[\int_\R \chi_{[0,\lambda)}(s) (V_* dx) = \int_\R \chi_{[0,\lambda)}(s) ((V_0)_* dx) \]
for every $\lambda > 0$, implying that $V_*dx$  agrees with an absolutely continuous measure on a family of sets which is generating for the Borel subsets of $(-\infty,\lambda_0)$, and in particular that the critical set of $V$ is Lebesgue-null.

\end{document}